\documentclass[11pt,leqno,twoside]{article}
\def\cvd{\hfill$\Box$}

\def\w{\dot{w}}

\def\int{\mathbb{Z}}
\def\C{\mathbb{C}}
\def\R{\mathbb{R}}
\def\z{\mathbb{Z}}
\def\Z{\mathbb{Z}}

\def\Q{\mathbb{Q}}

\def\Aut{\text{Aut}~}
\def\Diagr{\text{Diagr}}

\def\LAut{\rm{LAut}}
\def\AUT{\rm{AUT}}
\def\Lie{{\rm{Lie}}}
\def\Ad{{\rm{Ad}}~}

\def\O{{\cal O}}

\def\A{{\cal A}}

\def\gen{\mathfrak g}
\def\nen{\mathfrak n}
\def\uen{\mathfrak u}
\def\hen{{\mathfrak h}}

\def\len{\mathfrak l}
\def\men{\mathfrak m}

\def\K{{K}}

\def\slen{\mathfrak {sl}}

\def\proof{\noindent{\bf Proof. }}

\def\pf{\proof}

\def\rk{{\rm rk}}

\def\a{\alpha}
\def\f{\varphi}
\def\b{\beta}
\def\d{\delta}
\def\D{\Delta}

\def\g{\gamma}

\def\<#1{\langle #1\rangle}

\usepackage{times}
\usepackage{graphics}
\usepackage{amssymb}
\usepackage{amscd}
\usepackage{amsmath}
\usepackage{fancyhdr}

\input xy
\xyoption{all}

\title{Local automorphisms of finite dimensional\\ simple Lie algebras}

\newtheorem{theorem}{Theorem}[section]
\newtheorem{lemma}[theorem]{Lemma}
\newtheorem{corollary}[theorem]{Corollary}
\newtheorem{proposition}[theorem]{Proposition}
\newtheorem{definition}[theorem]{Definition}
\newtheorem{remark}[theorem]{Remark}

\newcounter{tigre}

\fancypagestyle{plain}{%
\fancyhf{}

}

\author{Mauro Costantini\\
Dipartimento di Matematica ``Tullio Levi-Civita''\\
Torre Archimede - via Trieste 63 - 35121 Padova - Italy\\
email: costantini@math.unipd.it }

\date{}
\begin{document}

\baselineskip=20pt
\maketitle
\begin{abstract}
  Let ${\mathfrak g}$ be a finite dimensional simple Lie algebra over an
  algebraically closed field $\K$ of characteristic $0$. A linear map $\f:{\mathfrak g}\to {\mathfrak g}$ is called a local automorphism if for every $x$ in ${\mathfrak g}$ there is an automorphism $\f_x$ of ${\mathfrak g}$ such that $\f(x)=\f_x(x)$. We prove that a linear map $\f:{\mathfrak g}\to {\mathfrak g}$ is local automorphism if and only if it is an automorphism or an anti-automorphism.
  
\medskip
\medskip
  
\medskip

\medskip

\end{abstract}
\section{Introduction}

\newcounter{equat}
\def\theequat{(\arabic{equat})}
\def\equat{\refstepcounter{equat}$$~}
\def\endequat{\leqno{\boldsymbol{(\arabic{equat})}}~$$}

\newcommand{\elem}[1]{\stackrel{#1}{\longto}}
\newcommand{\map}[1]{\stackrel{#1}{\to}}
\def\imp{\Rightarrow}
\def\Imp{\Longrightarrow}
\def\iff{\Leftrightarrow}
\def\Iff{\Longleftrightarrow}
\def\to{\rightarrow}
\def\longto{\longrightarrow}
\def\injto{\hookrightarrow}
\def\rtordu{\rightsquigarrow}

Mappings which are close to automorphisms and derivations of algebras have been
extensively investigated: in particular, since the 1990s (see \cite{Ka}, \cite{LS}, \cite{Se}), the description of local and $2$-local automorphisms (respectively, local and $2$-local derivations) of algebras 
has been deeply studied by many authors. 

Given an algebra $\A$ over a field $k$, a {\it local automorphism} (respectively, {\it local derivation}) of $\A$ is a $k$-linear map $\f:\A\to\A$ such that for each $a\in\A$ there exists an automorphism (respectively, a derivation)  $\f_a$ of $\A$ such that $\f(a)=\f_a(a)$.
A map $\f:\A\to\A$ (not $k$-linear in general) is called a {\it 2-local automorphism} (respectively, a {\it 2-local derivation}) if for every $x$, $y\in\A$, there exists an automorphism (respectively, a derivation) $\f_{x,y}$ of $\A$ such that $\f(x)=\f_{x,y}(x)$ and $\f(y)=\f_{x,y}(y)$.

In \cite{LS} the author proves that the automorphisms and the anti-automorphisms of the associative algebra $M_n(\C)$ of complex $n\times n$ matrices exhaust all its local automorphisms. On the other hand, it is proven in \cite{Cr} that a certain commutative subalgebra of $M_3(\C)$ has a local automorphism which is not an automorphism. 

Among other results (see the Introduction of \cite{AK4} for a detailed historical account), assuming the field $k$ is algebraically closed of characteristic zero, in \cite{AK1} the authors proved that every 2-local derivation of a finite dimensional semisimple Lie algebra is a derivation; in \cite{AK2} it is proved that every local derivation of a finite dimensional semisimple Lie algebra is a derivation. As far as automorphisms are concerned, in \cite{CW} the authors proved that if $\gen$ is a finite dimensional simple Lie algebra of type $A_\ell$ ($\ell\geq 1$), $D_\ell$ ($\ell\geq 4$), or $E_i$ ($i=6$, 7, 8), then every 2-local automorphism of $\gen$ is an automorphism. This result was extended to any finite dimensional semisimple Lie algebra in \cite{AK3}. On the other hand, for local automorphisms of simple Lie algebras it is only known that 
the automorphisms and the anti-automorphisms of the special linear algebra $\slen(n)$ exhaust all its local automorphisms (\cite[Theorem 2.3]{AK4}).

The main purpose of this paper it to extend this result to any finite dimensional simple Lie algebra: namely we prove that a $\K$-linear endomorphism of a finite dimensional simple Lie algebra $\gen$ over the algebraically closed field $\K$ of characteristic zero is a local automorphism if and only if it is an automorphism or an anti-automorphism of $\gen$.

Let $G$ be the connected component of the automorphism group of $\gen$: then $G$ is the adjoint simple algebraic group over $\K$ with the same Dynkin diagram as $\gen$. It is clear that every automorphism of $\gen$ is a local automorphism: we show that every anti-automorphism of $\gen$ is a local automorphism too. For this purpose we make use of the Bala-Carter theory for the classification of nilpotent elements in $\gen$.

To show that a local automorphism of $\gen$ is an automorphisms or an anti-automorphisms, we make use of the Tits' Building $\Delta(G)$ of $G$ (as definend in \cite[Chap. 5.3]{TITS}) and the classification theorem \cite[Theorem 5.8]{TITS} which in particular describes the automorphisms of $\Delta(G)$.

\setcounter{equation}{0}
\section{Preliminaries}\label{bepi}

Throughout the paper $\K$ is an algebraically closed field of characteristic zero. We denote by $\R$ the reals, by
$\Z$ the integers.

Let $A = (a_{ij})$ be a finite indecomposable
Cartan matrix of rank $n$. 
To $A$ there is associated a root system $\Phi$, a simple Lie 
algebra
${\mathfrak g}$ and a simple adjoint algebraic group $G$ over 
$\K$. We fix a maximal torus $T$ of $G$, and a Borel subgroup $B$
 containing $T$:  $B^-$ is the Borel subgroup opposite to $B$, $U$ (respectively $U^-$) is the unipotent radical of $B$ (respectively of $B^-$).  We denote by $\hen$, $\nen$, $\nen^-$ the Lie
algebra of $T$, $U$, $U^-$ respectively. Then  $\Phi$ is the set of roots relative
to $T$, and $B$ determines the set of positive roots  
$\Phi^+$, and the simple roots $\Delta=\{\alpha_1,\ldots,\alpha_n\}$.  The real space $E=\R \Phi$  is a Euclidean space, endowed with the scalar product
$(\a_i,\alpha_j) = d_ia_{ij}$. Here
$\{d_1,\ldots,d_n\}$  are
relatively prime positive integers such that if $D$ is the
diagonal matrix with entries $d_1,\ldots,d_n$, then $DA$ is
symmetric. For $\b=m_1\a_1+\cdots+m_n\a_n$, the height of $\b$ is $m_1+\cdots+m_n$. For $\a$, $\b\in\Phi$, we put ${\<{\b,\a}}=\frac{2(\b,\a)}{(\a,\a)}$.

We denote by $W$ the Weyl
group; $s_\a$ is the reflection associated to $\alpha\in\Phi$, we write for short $s_i$ for the simple reflection associated to $\alpha_i$, $w_0$ is the longest element of $W$. 
We put
$\Pi=\{1,\ldots,n\}$, $\vartheta$ the symmetry (called the opposite involution) of $\Pi$ induced by $-w_0$ and we fix a Chevalley basis $\{h_{i},  i\in \Pi; e_\a,\a\in\Phi\}$
 of $\frak g$ (see \cite[Chap. 4.2]{Carter1}). We put $h_\b=[e_\b,e_{-\b}]$ for $\b\in \Phi$ (hence $h_i=h_{\a_i}$ for $i\in \Pi$).

We use the notation $x_\a(\xi)$, for $\a\in \Phi$, $\xi\in \K$, as in \cite{Carter1}, \cite{springer}. For $\a\in \Phi$ we put $X_\a=\{x_\a(\xi)\mid \xi\in \K\}$, the root-subgroup corresponding to $\a$. We identify $W$ with $N/T$, where $N$ is the normalizer of $T$.

We choose the $x_\a$'s  so that, for all $\a\in \Phi$, $n_\a=x_\a(1)x_{-\a}(-1)x_\a(1)$
lies in $N$ and has image the reflection $s_\a$ in $W$. 
The family $(x_\a)_{\a\in\Phi}$ is called a {\it realization} of $\Phi$ in $G$.

Given an element $w\in W$ we shall denote
a representative of $w$ in $N$ by $\dot{w}$. We can, and shall, take $\dot w$ defined over $\Z$.

For algebraic groups we use the notation in \cite{Hu2}, \cite{Carter2}. In particular, 
for $J\subseteq \Pi$, $\Delta_J=\{\a_j\mid j\in J\}$, $\Phi_J$ is the corresponding root system, $W_J$ the Weyl group, $P_J$ the standard parabolic subgroup of $G$, $L_J=T\<{X_\a\mid \a\in \Phi_J}$ the standard Levi subgroup of $P_J$.
For $w\in W$ we have
$$
\dot wU^-\dot w^{-1}\cap U=\prod_{\substack{\a>0\\w^{-1}\a<0}}X_\a
$$

If $x$ is an element of $\gen$, $C_G(x)$ is the centralizer of $x$ in $G$. 

We denote by $GL(\gen)$ the group of automorphisms of $\gen$ as a $\K$-vector space.
The group $\AUT(\gen)$ of automorphisms of $\gen$ as a Lie algebra is completely described in \cite[Chap. IX]{Jac}, \cite[16.5]{Hu1}.

We denote by $\cal NB(\gen)$ the set of the nilradicals of Borel subalgebras of $\gen$. This is a unique orbit under $G$: if $\nen_1\in \cal NB(\gen)$ then, by the Bruhat decompositon in $G$, there exists a unique $w\in W$ and a unique $u\in \dot wU^-\dot w^{-1}\cap U$ such that
${\nen}_1=\Ad u\dot w.\nen$.

\setcounter{equation}{0}
\section{The main result}\label{main}

We recall that a parabolic subgroup $P$ is called {\it distinguished} if $\dim P/R_uP=\dim R_uP/(R_uP)'$. Here  $R_uP$ is the unipotent radical of $P$ and $(R_uP)'$ is the derived subgroup of $R_uP$ (see \cite[p. 167]{Carter2}). Two parabolic subgroups are said to be {\it opposite} if their intersection is a common Levi subgroup (see \cite[14.20]{Bo}). If $P$ is a parabolic subgroup and $L$ is a Levi subgroup of $P$, then there exists a unique parabolic subgroup opposite to $P$ containing $L$.  Any two opposite parabolic subgroups of $P$ are conjugate by a unique element of $R_uP$ (\cite[Proposition 14.21]{Bo}). 

\begin{lemma}\label{opposite} Let $P$ be a distinguished parabolic subgroup of a semisimple algebraic group $R$ and let $P^{\rm op}$ be an opposite parabolic subgroup of $P$. Then $P$ and $P^{\rm op}$ are conjugate in $R$.
\end{lemma}
\pf It is enough to assume $R$ simple, $P=P_J=\<{B, X_{-\a_i}\mid i\in J}$, $P^{\rm op}=\<{B^-, X_{\a_i}\mid i\in J}$ for a certain $J\subseteq \Pi$. 
If $w_0 = -1$, then $P^{\rm op} = \dot w_0 P \dot w_0^{-1}$. We are left with the cases where $R$ is of type $A_n$, $n\geq 2$, $D_n$ with $n\geq 5$, $n$ odd, $E_6$. From the tables in \cite{Carter2}, p. 174 - 176, one checks that again $P^{\rm op} = \dot w_0 P \dot w_0^{-1}$, since in each case the diagram of $P$ is invariant under the opposite involution $\vartheta$ of the Dynkin diagram.
\cvd

\begin{theorem}\label{inversa} The anti-automorphism $-i_\gen:\gen\to\gen$, $x\mapsto -x$ is a local automorphism of $\gen$.
\end{theorem}
\pf 
Let $x\in\gen$. We have to show that there exists $\a\in\AUT(\gen)$ such that $\a(x)=-x$. 
Let $\O$ be the $G$-orbit of $x$: it is enough to show that this holds for a certain $y\in \O$. In fact, if $x=\Ad g.y$ and $\b(y)=-y$ for certain $g\in G$, $\b\in\AUT(\gen)$, then $\a(x)=-x$, where $\a$ is the automorphism of $\gen$ given by $\a=(\Ad g)\b(\Ad g)^{-1}$.
 
 Let $x=s+e$ be the Jordan-Chevalley decomposition of $x$, i.e.  $s$ is semisimple, $e$ is nilpotent, with $[s,e]=0$. Let $H=C_G(s)$. This is a Levi subgroup of $G$ and, up to conjugacy in $G$, we may assume that $H$ is the standard Levi subgroup $L_J$ of $G$. Moreover the centralizer of $s$ in $\gen$ is the Lie algebra $\len_J$ of $L_J$, $e$ lies in $\len_J$, and $s$ lies in the center $Z(\len_J)\subseteq \hen$. Let $\men$ be a minimal Levi subalgebra of $\len_J$ containing $e$. Let $M$ be the Levi subgroup of $H$ such that $\men=\Lie(M)$, and let $M'$ be the semisimple part of $M$ and $\men'=\Lie(M')$. Then $e$ lies in $\men'$ and $e$ is distinguished in $\men'$. There exists a distinguished parabolic subgroup $P_{M'}$ of $M'$ such that $e$ lies in the dense orbit of $P_{M'}$ on the Lie algebra $\uen_{P_{M'}}$ of its unipotent radical. Up to conjugation by an element of $H$, we may assume that $P_{M'}=\<{T_1, X_\a, X_{-\a}, X_\delta\mid \a\in \Psi_1, \delta\in\Psi_2 }$ for $T_1$ a certain subtorus of $T$ and $\Psi_1$, $\Psi_2$ certain disjoint subsets of $\Phi^+$.
 
 Now we consider an automomorphism $\imath$ of $G$ satisfying
 $$
\imath(t)=t^{-1}
 \quad \text{for every}\  t\in T\quad,\quad
\imath(X_\a)=X_{-\a}
  \quad \text{for every}\  \a\in \Phi
 $$
 \cite[proof of Corollary 1.16, p. 189]{Jan}.
 Then the differential $d\imath$ is an automorphism of $\gen$ satisfying $d\imath(h)=-h$ for every $h\in \hen$, in particular $d\imath(s)=-s$. It is enough, by \cite[Lemma 2.2.1]{ST}, to show that $d\imath(e)$ and $-e$ are conjugate by an element of $H$. But $\imath(P_{M'})=\<{T_1, X_\a, X_{-\a}, X_{-\delta}\mid \a\in \Psi_1, \delta\in\Psi_2 }$ is opposite to $P_{M'}$ (since $P_{M'}\cap \imath(P_{M'})=\<{T_1, X_\a, X_{-\a}
 \mid \a\in \Psi_1 }$, a Levi subgroup of $M'$). Since a parabolic subgroup has a unique dense orbit on the Lie algebra of its unipotent radical, and clearly $-e$ lies
in the dense orbit of $P_{M'}$ on $\uen_{P_{M'}}$, and $d\imath(N)$ lies in the dense orbit of $\imath(P_{M'})$ on $d\imath(\uen_{P_{M'}})$, its is enough to show that 
$P_{M'}$ and $\imath(P_{M'})$ are conjugate in $H$. From Lemma \ref{opposite} it follows that $P_{M'}$ and $\imath(P_{M'})$ are already conjugate in $M'$, and we are done.\cvd

 \medskip
 We denote by $\AUT^\ast\!(\gen)$ the group of automorphisms of the $\K$-vector space $\gen$ which are either automorphisms or anti-automorphisms of the Lie algebra $\gen$.
 Then $\AUT^\ast\!(\gen)=\AUT(\gen)\rtimes \<{-i_\gen}$. We observe that if $\f$ is a local automorphism of $\gen$, then $\f$ is invertible and its inverse is a local automorphism. It is also clear that the composite of local automorphisms is a local automorphism, therefore the set $\LAut(\gen)$ of local automorphisms of $\gen$ is 
 a subgroup of $GL(\gen)$.
  By Theorem \ref{inversa} we have
  
  \begin{corollary}\label{anti} Every anti-automorphism of $\gen$ is a local automorphism, i.e.  $\AUT^\ast\!(\gen)\leq \LAut(\gen)$.
  \phantom{abcde}\hfill
  \cvd
  \end{corollary}
   We shall prove that  $\LAut(\gen)= \AUT^\ast\!(\gen)$.

 \begin{lemma}\label{invarianza} Let $\f$ be in $\LAut(\gen)$. Then $\f$ leaves invariant the set $\cal N$ of nilpotent elements and the set $\cal S$ of semisimple elements of $\gen$.
 \end{lemma}
\pf Let $x\in\gen$. There exists $\f_x\in\AUT(\gen)$ such that $\f_x(x)=\f(x)$. Since automorphisms map nilpotent (respectively, semisimple) elements to nilpotent (respectively, semisimple) elements, it follows that $\f(\cal N)\subseteq \cal N$ and $\f(\cal S)\subseteq \cal S$. Since $\f^{-1}$ is also a local automorphism, we conclude that 
$\f(\cal N)= \cal N$ and $\f(\cal S)= \cal S$.\cvd

\medskip

A classical theorem of Gerstenhaber \cite{Ge} states that any vector space consisting of nilpotent $n\times n$ matrices has dimension at most $\frac12 n(n-1)$, and that any such space attaining this maximal
possible dimension is conjugate to the space of upper triangular matrices. In \cite{DKK} the authors generalized this result to the Lie algebra of any reductive algebraic group over any algebraically closed field, under certain conditions in case the characteristic of the field is 2 or 3. We restate this generalization for our purposes. For short we say that a subspace $V$ of $\gen$ is nilpotent, if $V$ consists of nilpotent elements.

\begin{theorem}\label{sottoalgebre}(\cite[Theorem 1]{DKK})
Let $V$ be a nilpotent subspace of a finite dimensional semisimple Lie algebra $\gen$ over $\K$. Then $\dim V  \leq \frac12 (\dim \gen-\rk\,  \gen)$ and, if equality holds, $V$ is the nilradical
of a Borel subalgebra of $\gen$.
\end{theorem}

In particular the nilpotent subspaces of maximal dimension are the maximal nilpotent subalgebras $\gen$: they constitute the set $\cal NB(\gen)$ defined in the Preliminaries.

 \begin{proposition}\label{orbita} Let $\f$ be in $\LAut(\gen)$. Then $\f$ induces a permutation of the set $\cal NB(\gen)$.
\end{proposition}
\pf Let $V$ be any nilpotent subspace of $\gen$. By Lemma \ref{invarianza} $\f(V)$ and $\f^{-1}(V)$ are nilpotent subspaces of $\gen$. Therefore $\f$ induces a permutation $V\mapsto \f(V)$ of the set of all nilpotent subspaces of $\gen$. In particular $\f$ induces a permutation of $\cal NB(\gen)$.\cvd

\medskip

We introduce the canonical Tits' Building $\Delta(G)$ associated to $G$.
 
 \begin{definition}\label{building}\cite[Chap. 5.3]{TITS} The building $\Delta(G)$ of $G$ is the set of all parabolic subgroups of $G$, partially ordered by reverse of inclusion.
 \end{definition}
 
 The maximal elements of $\D(G)$ (called {\it chambers}) are the Borel subgroups of $G$. The set of Borel subgroups of $G$ is in canonical bijection with the set of 
 Lie algebras of Borel subgroups of $G$ (i.e. the Borel subalgebras of $\gen$, \cite[14.25]{Bo}), and this set is in canonical bijection with the set $\cal NB(\gen)$. By Proposition \ref{orbita}, a local automorphism $\f$ of $\gen$ induces a permutation of $\cal NB(\gen)$, and therefore a permutation $\rho_\f$ of the set of chambers of $\Delta(G)$. Let $B_1$, $B_2$ be adjacent chambers: this means that the codimension (as algebraic varieties) of $B_1\cap B_2$ in $B_1$ (and $B_2$) is 1. Since $B_1\cap B_2$ always contains a maximal torus of $G$, this is equivalent to the condition that the codimension (as $k$-vector spaces) of $\nen_1\cap \nen_2$ in $\nen_1$ (and $\nen_2$) is 1, where $\nen_i$ is the nilradical of the Lie algebra of $B_i$ for $i=1$, $2$.
  
\begin{proposition}\label{estensione} Let $\f$ be in $\LAut(\gen)$. Then $\rho_\f$ can be (uniquely) extended to an automorphism of $\D(G)$.
\end{proposition}
\pf By the previous discussion, this follows from \cite[Theorem 3.21, Corollary 3.26]{TITS}.\cvd

\noindent
We shall still denote by $\rho_\f$ the automorphism of $\D(G)$ induced by $\f$.

\medskip

A {\it symmetry} of the Dynkin diagram of $G$ is a permutation $\d$ of the nodes of the diagram such that ${\<{\a_{\d(i)},\a_{\d(j)}}}={\<{\a_i,\a_j}}$ for all $i$, $j\in \Pi$ (\cite[p. 277]{Jac}. Note that in \cite[p. 200]{Carter1} the definition is different, in order to deal also with fields of characteristic 2 or 3). We denote the group of symmetries of the Dynkin diagram by $\Diagr$.

\begin{definition}\label{diagramma}  Let $\d$ be a symmetry of the Dynkin diagram of $\gen$ . We denote by $d_\d$ both the isometry of $E$ and the {\it graph automorphism} of $\gen$ defined respectively by
$$
d_\d({\a_i})={\a_{\d (i)}}\ \quad  \text{ for every}\  i\in \Pi
$$
$$
d_\d(e_{\a_i})=e_{\a_{\d (i)}}\ , \ d_\d(e_{-\a_i})=e_{-\a_{\d (i)}}\ , \ d_\d(h_{\a_i})=h_{\a_{\d (i)}}\quad  \text{ for every}\  i\in \Pi
$$
\end{definition}
 
 \begin{proposition}\label{scalari} Let $\f=c\, i_\gen$, for a certain $c\in \K^\ast$. Then $\f\in \LAut(\gen)$ if and only if $c=\pm 1$.
\end{proposition}
\pf We only need to show that if $\f=c\, i_\gen$ is a local automorphism, then $c=\pm 1$. By \cite[Proposition 6.4.2]{Carter1} we have
$$
\Ad n_\a.h_\b=h_{s_\a(\b)}
$$
for every $\a$, $\b\in\Phi$, so that 
$$
\Ad \dot w.h_\b=h_{w(\b)}
$$
for every $w\in W$, $\b\in\Phi$. Now fix any $\a\in \Phi$, $h\in\hen$. There exists $g\in G$, $\delta\in \Diagr$ such that 
$$
c\,h_\a=\f(h_\a)=d_\delta \Ad g.h_\a
$$
Hence $\Ad g.h_\a = c\, d^{-1}_\d h_\a\in \hen$, which means that the elements $h_\a$ and $c\, d^{-1}_\d h_\a$ of $\hen$ are conjugate under $G$, and therefore they are conjugate under $W$, i.e. there exists $w\in W$ such that $\Ad g.h_\a = \Ad \dot w. h_\a = h_{w(\a)}$. Hence
$c\, d^{-1}_\d h_\a=h_{w(\a)}$, $c\, h_\a=d_\d h_{w(\a)}=h_{\d w(\a)}=h_\b$, for $\b=\d w(\a)\in\Phi$. It follows that $\b=\pm \a$, i.e. $c=\pm1$.\cvd

\medskip
A {\it semilinear isomorphism} between two Lie algebras is a bijective semilinear mapping of the underlying vector spaces which respects Lie multiplication. 

\begin{definition} Let $f \in \Aut \K$. We denote by $a_f$ both the {\it field automorphism} of $G$ (as an abstract group) and the $f$-semilinear automorphism of $\gen$ defined respectively by
$$
a_f(x_\a(k))=x_\a(f(k))\quad  \text{ for every}\  \a\in \Phi, k\in \K
$$
$$
a_f(k e_\a)=f(k)e_\a\quad  \text{ for every}\  \a\in \Phi, k\in \K
$$
\end{definition}

\begin{remark}{\rm 
Note that we also have $a_f(k h_\a)=f(k)h_\a$ for every $\a\in \Phi$, $k\in \K$, 
since $h_\a=[e_\a,e_{-\a}]$ for every $\a\in\Phi$.
 Moreover, for every $g\in G$, $x\in \gen$ we have $a_f(\Ad g.x)=\Ad (a_f(g)).a_f(x)$.
}
\end{remark}

\begin{proposition}\label{campo} 
Let $\f \in GL(\gen)$ and $f \in \Aut \K$ be such that $\f(X)=a_f(X)$ for every $X\in \cal NB(\gen)$. Then $f=i_\K$ and there is $c\in \K^\ast$ such that $\f=c\, i_\gen$.
\end{proposition}
\pf 
We have $a_f(\nen)=\nen$ and 
$a_f(\nen^-)=\nen^-$. It follows that
$$
a_f(\Ad x_\a(k)\dot w.\nen)=\Ad x_\a(f(k))\dot w.\nen\quad,\quad
a_f(\Ad x_\a(k)\dot w.\nen^-)=\Ad x_\a(f(k))\dot w.\nen^-
$$
for every $\a\in \Phi$, $k\in \K$, since we fixed the representatives $\dot w$ over $\z$, and therefore $a_f(\dot w)=\dot w$ for every $w\in W$.

We shall repeatedly use the fact that if $\nen_1$, $\nen_2\in \cal NB(\gen)$ are such that $\nen_1\cap \nen_2=\<v$ with $v\not=0$, then
$$
\<{\f(v)}=\f(\nen_1)\cap \f(\nen_2)=a_f(\nen_1)\cap a_f(\nen_2)=\<{a_f(v)}
$$

For every $i\in \Pi$ we have 
$$
\Ad \dot s_i.\nen^-\cap \nen = \<{e_{\a_i}}\quad,\quad
\Ad \dot s_i.\nen\cap \nen^- = \<{e_{-\a_i}}
$$
hence
$$
 \<{\f(e_{\a_i})}= \<{a_f(e_{\a_i})}= \<{e_{\a_i}}\quad,\quad
  \<{\f(e_{-\a_i})}= \<{a_f(e_{-\a_i})}= \<{e_{-\a_i}}
$$
Let $\a\in \Phi$. There exists $w\in W$, $i\in \Pi$ such that $w(\a_i)=\a$. Then 
$$
\<{e_{\a}}=\Ad \dot w.\<{e_{\a_i}}=\Ad \dot w\dot s_i.\nen^-\cap \Ad \dot w.\nen 
$$
so that
$$
 \<{\f(e_{\a})}= \<{a_f(e_{\a})}= \<{e_{\a}}
$$
Hence, for every $\a\in \Phi$ there exists $c_\a\in \K^\ast$ such that $\f(e_\a)=c_\a e_\a$. 

By \cite[p. 64]{Carter1}, for every $\a\in \Phi$, $k\in \K$ we have
$$
\Ad x_{\a}(k).e_\a=e_\a\quad,\quad
\Ad x_{\a}(k).e_{-\a}=e_{-\a}+kh_\a-k^2e_\a
$$
Let us fix $i$ in $\Pi$. From
$
\Ad \dot s_i.\nen\cap \nen^- = \<{e_{-\a_i}}
$
we get 
$$
\Ad x_{\a_i}(k).\Ad \dot s_i.\nen\cap \Ad x_{\a_i}(k).\nen^-=\<{\Ad x_{\a_i}(k).e_{-\a_i}}=\<{e_{-\a_i}+kh_{\a_i}-k^2e_{\a_i}}
$$
so that 
\begin{equation}\label{generico}
\begin{split}
\<{\f(e_{-\a_i}+kh_{\a_i}-k^2e_{\a_i})}&=\<{a_f(e_{-\a_i}+kh_{\a_i}-k^2e_{\a_i})}=\\
&=\<{e_{-\a_i}+f(k)h_{\a_i}-f(k)^2e_{\a_i}}
\end{split}
\end{equation}
In particular, for $k=1$ we get
$$
\<{\f(e_{-\a_i}+h_{\a_i}-e_{\a_i})}
=
\<{e_{-\a_i}+h_{\a_i}-e_{\a_i}}
$$
hence $\f(h_{\a_i})=d_i h_{\a_i}+x_i e_{\a_i}+y_i e_{-\a_i}$ for certain $d_i$, $x_i$, $y_i\in \K$, $i=1,\ldots,n$. From (\ref{generico}), for every $k\in \K$ there exits $p_k\in \K^\ast$ such that
\begin{equation}\label{combinazione}
c_{-\a_i}e_{-\a_i}+k(d_i h_{\a_i}+x_i e_{\a_i}+y_i e_{-\a_i})-k^2c_{\a_i}e_{\a_i}=p_k(e_{-\a_i}+f(k)h_{\a_i}-f(k)^2e_{\a_i})
\end{equation}
hence
$
k\, d_i =p_kf(k)
$
for every $k\in \K$ and in particular, for $k=1$, $d_i=p_1$. But then $p_k=\frac{k}{f(k)}p_1$ for every $k\in \K^\ast$, so that $p_k=p_1$ for every $k$ in the prime field $\Q$ of $\K$, $k\not=0$. From (\ref{combinazione}) we obtain $c_{-\a_i}e_{-\a_i}+k y_i e_{-\a_i}=p_1e_{-\a_i}$ 
and
$
k x_i e_{\a_i}-k^2c_{\a_i}e_{\a_i}=-p_1 k^2e_{\a_i}$
for every $k\in\Q^\ast$, so that $y_i=0$, $c_{-\a_i}=p_1$, $x_i=0$ and $c_{\a_i}=p_1$. We have proved that
\begin{equation}\label{semplici}
\f(h_{\a_i})=c_{\a_i} h_{\a_i}\ , \ c_{-\a_i}=c_{\a_i}
\end{equation}
Moreover, from (\ref{combinazione}) it follows that
$$
c_{\a_i}e_{-\a_i}+kc_{\a_i} h_{\a_i}-k^2c_{\a_i}e_{\a_i}=p_k(e_{-\a_i}+f(k)h_{\a_i}-f(k)^2e_{\a_i})
$$
for every $k\in \K$, hence $p_k=c_{\a_i}$ and $f(k)=k$ for every $k\in \K$, i.e. $f=i_\K$.

So far we have proved that $f=i_\K$, and that for every $i=1,\ldots,n$ we have 
$\f(e_{\a_i})=c_{\a_i} e_{\a_i}$, 
$\f(e_{-\a_i})=c_{\a_i} e_{-\a_i}$ and
$\f(h_{i})=c_{\a_i} h_{i}$.
Our aim is to show that $c_\a=c_\b$ for every $\a$, $\b\in \Phi$. We prove that $c_\a=c_\b$ for every $\a$, $\b\in \Phi^+$. With a similar procedure it will follow that  $c_\a=c_\b$ for every $\a$, $\b\in \Phi^-$, so that  $c_\a=c_\b$ for every $\a$, $\b\in \Phi$ by (\ref{semplici}).

By \cite[p. 64]{Carter1}, for linearly independent roots $\a$, $\b$ we have
$$
\Ad x_\a(t).e_\b=\sum_{r=0}^q M_{\a,\b,r}\,t^r \,e_{r\a+\b}
$$
where $M_{\a,\b,0}=1$,  $M_{\a,\b,r}=\pm {p+r\choose r}$ for $r\geq 1$, $-p\a+\b,\ldots,\b,\ldots,q\a+\b$ is the $\a$-chain through $\b$ with $p$ and $q$ non negative integers. In particular, for $t=1$ we get
\begin{equation}\label{indipendenti}
\Ad x_\a(1).e_\b=\sum_{r=0}^q M_{\a,\b,r} \,e_{r\a+\b}
\end{equation}
We begin by showing that $c_{\a_i}=c_{\a_j}$ for every $i$, $j\in \Pi$. Assume $\a_i+\a_j\in\Phi$. Then
$$
\Ad x_{\a_i}(1).e_{\a_j}=\sum_{r=0}^q M_{\a_i,\a_j,r} \,e_{r\a_i+\a_j}
$$
with $q\geq 1$.
From $\Ad \dot s_j.\nen^-\cap \nen = \<{e_{\a_j}}$ we get
$$
\<{\Ad x_{\a_i}(1).e_{\a_j}}=\Ad x_{\a_i}(1)\Ad \dot s_j.\nen^-\cap \Ad x_{\a_i}(1).\nen
$$
so that
$$
\<{\f(\Ad x_{\a_i}(1).e_{\a_j})}=\<{a_f(\Ad x_{\a_i}(1).e_{\a_j})}=\<{\Ad x_{\a_i}(1).e_{\a_j}}
$$
There exists $c\in \K^\ast$ such that
$$
\f(\sum_{r=0}^q M_{\a_i,\a_j,r} \,e_{r\a_i+\a_j})=c\,(\sum_{r=0}^q M_{\a_i,\a_j,r} \,e_{r\a_i+\a_j})
$$
Since $M_{\a_i,\a_j,r}\not=0$ for every $r=0,\ldots,q$, we get 
$$
c_{r\a_i+\a_j}=c
$$
for every $r=0,\ldots,q$, and in particular $c_{\a_j}=c$,  $c_{\a_i+\a_j}=c$, so that $c_{\a_j}=c_{\a_i+\a_j}=c$.
Similarly, by considering $\Ad x_{\a_j}(1).e_{\a_i}$, we obtain $c_{\a_i}=c_{\a_j+\a_i}$: hence $c_{\a_i}=c_{\a_j}=c$.
Since the Dynkin diagram is connected, we get $c_{\a_i}=c_{\a_j}=c$ for every $i$, $j\in \Pi$ (incidentally, the previous argument shows that $c_{\a}=c_{\b}=c$ for positive roots $\a$, $\b$ of height at most 2).

Assume that $\b$ is a positive root of height $m$ with $m\geq 2$. Then we may write $\b=\g+\a_i$, for a certain $\g\in\Phi^+$ (of height $m-1)$ and a certain $i\in\Pi$. Then
$$
\Ad x_{\g}(1).e_{\a_i}=\sum_{r=0}^q M_{\g,\a_i,r} \,e_{r\g+\a_i}
$$
with $q\geq 1$.
From $\Ad \dot s_i.\nen^-\cap \nen = \<{e_{\a_i}}$ we get
$$
\<{\Ad x_{\g}(1).e_{\a_i}}=\Ad x_{\g}(1)\Ad \dot s_i.\nen^-\cap \Ad x_{\g}(1).\nen 
$$
so that
$$
\<{\f(\Ad x_{\g}(1).e_{\a_i})}=\<{a_f(\Ad x_{\g}(1).e_{\a_i})}=\<{\Ad x_{\g}(1).e_{\a_i}}
$$
There exists $d\in \K^\ast$ such that
$$
\f(\sum_{r=0}^q M_{\g,\a_i,r} \,e_{r\g+\a_i})=d\,(\sum_{r=0}^q M_{\g,\a_i,r} \,e_{r\g+\a_i})
$$
Since $M_{\g,\a_i,r}\not=0$ for every $r=0,\ldots,q$, we get 
$$
c_{r\g+\a_i}=d
$$
for every $r=0,\ldots,q$, and in particular $c_{\a_i}=d$,  $c_\b=c_{\g+\a_i}=d$, so that $c_{\b}=c_{\a_i}=c$.
We have therefore proved that $c_\a=c$ for every $\a\in\Phi^+$. Similarly one can prove that $c_\a=c'$ for every $\a\in\Phi^-$ for a certain $c'\in\K^\ast$. Since 
by (\ref{semplici}) we have $c_{-\a_i}=c_{\a_i}$ we get $c'=c$, i.e. $c_\a=c$ for every $\a\in\Phi$. But we also have $\f(h_i)=c\,h_i$ for every $i\in\Pi$, we conclude that $\f=c\, i_\gen$.\cvd

\begin{theorem}\label{finale}
Let $\gen$ be a finite dimensional simple Lie algebra over the algebraically closed field $\K$ of characteristic zero. Then a linear map $\f:{\mathfrak g}\to {\mathfrak g}$ is local automorphism if and only if it is an automorphism or an anti-automorphism, i.e. $\LAut(\gen) =  \AUT^\ast\!(\gen)$.
\end{theorem}
\pf The case when $\gen$ is of type $A_n$, $n\geq 1$, is dealt with in \cite{AK4}. For completeness, here we give a proof also for this case. By Corollary \ref{anti} we have $\AUT^\ast\!(\gen)\leq \LAut(\gen)$.
Let $\f$ be a local automorphism of $\gen$. We show that there exists an automorphism $\b$ of $\gen$ and $c\in \K^\ast$ such that $\b^{-1}\f=c\, i_\gen$.

Assume first that $\gen$ has rank 1, i.e. $\gen=\slen(2)$. Then the result follows from the main theorem in \cite{BPW} (see Remark on page 45). So assume $\rk\, \gen \geq 2$. By Proposition \ref{estensione}, $\f$ induces an automorphism $\rho_\f$ of the building $\Delta(G)$ of $G$. By the structure theorem on isomorphisms of buildings (\cite[Theorem 5.8]{TITS}), there exists an automorphism $\a$ of $G$ (as an algebraic group) and a field automorphism $a_f$ of $G$ such that $\rho_\f(P)=\a\, a_f(P)$ for every parabolic subgroup $P$ of $G$. It follows that, for $\b=d\a$, the differential of $\a$, we get
$$
\b^{-1}\f(X)=a_f(X)
$$
for every $X$ in $\cal NB(\gen)$. By Proposition \ref{campo}, $\b^{-1}\f = c\, i_\gen$ for a certain $c\in\K^\ast$.

Finally, from Proposition \ref{scalari}, we get $c=\pm 1$, and $\f=\pm \b\in \AUT^\ast\!(\gen)$.
\cvd

\begin{remark}{\rm 
From the structure of the automorphism group of $\gen$, it follows that any $\f\in \LAut(\gen)$ is of the form $\f=\pm d_\d(\Ad g)$ for a unique $g\in G$ and a unique graph automorphism $d_\d$.
}
\end{remark}


\begin{thebibliography}{1}
     
\bibitem{AK1}{\sc S.\ Ayupov, K.\ Kudaybergenov, I.\ Rakhimov,}
\newblock{\em 2-local derivations on finite-dimensional {L}ie algebras,}
\newblock Linear Algebra Appl. 474, 1--11 (2015).

\bibitem{AK2}{\sc S.\ Ayupov, K.\ Kudaybergenov,}
\newblock{\em Local derivations on finite dimensional Lie algebras,}
\newblock Linear Algebra Appl. 493, 381--398 (2016).

\bibitem{AK3}{\sc S.\ Ayupov, K.\ Kudaybergenov,}
\newblock{\em 2-local automorphisms on finite-dimensional {L}ie algebras,}
\newblock Linear Algebra Appl. 507, 121--131 (2016).

\bibitem{AK4}{\sc S.\ Ayupov, K.\ Kudaybergenov,}
\newblock{\em Local automorphisms on finite-dimensional Lie and Leibniz algebras,}
\newblock preprint arXiv:1803.03142v2.

\bibitem{Bo}
{\sc A.\ Borel,}
\newblock {\em Linear Algebraic Groups,}
\newblock Second enlarged edition, Springer-Verlag, New York (1991).

\bibitem{BPW}
{\sc P.\ Botta, S.\ Pierce, W.\ Watkins,}
\newblock{\em Linear transformations that preserve the nilpotent matrices,}
\newblock  Pacific J. Math. 104(1), 39--46 (1983)

\bibitem{Carter1}
{\sc R.\ W.\ Carter,}
\newblock {\em Simple Groups of Lie Type,}
\newblock John Wiley (1989).

\bibitem{Carter2}
{\sc R.\ W.\ Carter,}
\newblock {\em Finite Groups of Lie Type,} 
\newblock John Wiley  (1985).

\bibitem{CW}
{\sc Z.\ Chen, D.\ Wang,}
\newblock{\em 2-Local automorphisms of finite-dimensional simple Lie algebras,}
\newblock  Linear Algebra Appl. 486, 335--344 (2015).

\bibitem{Cr}
{\sc R.\ Crist,}
\newblock{\em Local Automorphisms,}
\newblock  Proc. Amer. Math. Soc. 128, 1409--1414 (2000).

\bibitem{DKK}
{\sc J. \ Draisma, H.\ Kraft, J.\ Kuttler,}
\newblock{\em Nilpotent subspaces of maximal dimension in semi-simple {L}ie
              algebras,}
\newblock Compos. Math.\ 142(2), 464--476 (2006).


\bibitem{Ge}
{\sc M.\ Gerstenhaber,}
\newblock{\em On nilalgebras and linear varieties of nilpotent matrices, I,}
\newblock Amer. J. Math. 80, 614--622 (1958).

\bibitem{Hu1}
{\sc J.E.\ Humphreys,}
\newblock {\em Introduction to Lie algebras and representation theory,}
\newblock Third printing, Graduate Texts in Mathematics, No. 9, Springer-Verlag, New York-Heidelberg (1980).

\bibitem{Hu2}
{\sc J.E.\ Humphreys,}
\newblock {\em Linear Algebraic Groups,}
\newblock Third printing, Graduate Texts in Mathematics, No. 21, Springer-Verlag, New York-Heidelberg (1987).

\bibitem{Jac}
{\sc N.\ Jacobson,}
\newblock {\em Lie Algebras,}
\newblock Republication of the 1962 original, Dover Publications, Inc., New York (1979).

\bibitem{Jan}
{\sc J.C.\ Jantzen,}
\newblock {\em Representations of algebraic groups,}
\newblock Pure and Applied Mathematics 131, Academic Press, Inc., Boston, MA (1987).

\bibitem{Ka}
{\sc R.V.\  Kadison,}
\newblock{\em Local derivations,}
\newblock J. Algebra 130, 494--509 (1990).

\bibitem{LS}
{\sc D.R.\ Larson, A.R.\ Sourour,}
\newblock{\em Local Derivations and Local Automorphisms of B(H),}
\newblock {In: ``Operator theory: operator algebras and applications, {P}art 2''
              ({D}urham, {NH}, 1988)},
Proc. Sympos. Pure Math. 51, 187--194 (1990).

\bibitem{Se}
{\sc P.\ \v{S}emrl,}
\newblock{\em Local automorphisms and derivations on $B(H)$,}
\newblock Proc. Amer. Math. Soc. 125, 2677--2680 (1997).

\bibitem{ST}
{\sc A.\ Singh, M.\ Thakur,}
\newblock{\em Reality properties of conjugacy classes in algebraic groups,}
\newblock Israel J. Math. 165, 1--27 (2008).

\bibitem{springer}
{\sc T.A. Springer,}
\newblock{\em Linear Algebraic Groups,} Second Edition,
\newblock Progress in Mathematics 9, Birkh{\"a}user (1998).

\bibitem{TITS}
{\sc J.\ Tits,}
\newblock{\em Buildings of spherical type and finite {BN}-pairs,}
\newblock Lecture Notes in Mathematics, Vol. 386, Springer-Verlag, Berlin-New York (1974).


\end{thebibliography}
\end{document}